\documentclass{aiaa-tc}

\usepackage[sort&compress,numbers,square]{natbib}

\usepackage{setspace}
\usepackage{amsmath}
\usepackage{amsfonts}
\usepackage{amssymb,latexsym}
\usepackage[psamsfonts]{eucal}
\usepackage{amsthm}

\usepackage{xcolor}

\usepackage{BernsteinStyle}

\begin{document}



\title{\Large\bf{Euler's Equation via Lagrangian Dynamics\\ with Generalized Coordinates
}}

\author{Dennis S. Bernstein, University of Michigan, Ann Arbor, MI\footnote{Professor, Aerospace Engineering Department, Corresponding Author, \texttt{dsbaero@umich.edu}},\\ 
Ankit Goel, University of Maryland, Baltimore County, Baltimore, MD\footnote{Assistant Professor, Mechanical Engineering Department, \texttt{ankggoel@umbc.edu}}, \\
Omran Kouba, Higher Institute for Applied Sciences and Technology, Damascus, Syria\footnote{Professor, Department of Mathematics, \texttt{omran\_kouba@hiast.edu.sy}}}
\maketitle 


%
%

\begin{abstract}
Euler's equation relates the change in angular momentum of a rigid body to the applied torque.  This paper fills a gap in the literature by using Lagrangian dynamics to derive Euler's equation in terms of generalized coordinates.  This is done by parameterizing the angular velocity vector in terms of 3-2-1 and 3-1-3 Euler angles as well as Euler parameters, that is, unit quaternions.

\end{abstract}
 
\section{\large Introduction}

The rotational dynamics of a rigid spacecraft are modeled by Euler's equation \cite[p. 59]{hughes}, which relates the rate of change of the spacecraft angular momentum to the net torque.
Let $\omega\in\BBR^3$ denote the angular velocity of the spacecraft relative to an inertial frame, let $J\in\BBR^{3\times3}$ denote the inertia matrix of the spacecraft relative to its center of mass, and let $\tau$ denote the net torque applied to the spacecraft.  All of these quantities are expressed in the body frame.
Applying Newton-Euler dynamics yields Euler's equation  
\begin{equation}
    J\dot\omega + \omega\times J\omega = \tau. \label{ee}
\end{equation}

An alternative approach to obtaining the dynamics of a mechanical system is to apply Hamilton's principle in the form of Lagrangian dynamics given by
\begin{align}
    \rmd_t \p_{\dot q} T - \p_q T = Q, \label{ld}
\end{align}
where $T$ is the kinetic energy of the system, $q$ is the vector of generalized coordinates, and $Q$ is the vector of generalized forces arising from all external and dissipative forces and torques, including those arising from potential energy.
Here, $\rmd_t$ denotes total time derivative, and $\p_{\dot q}$ and $\p_q$ denote the partial derivatives with respect to $\dot q$ and $q,$ respectively.

For a mechanical system consisting of multiple rigid bodies, \eqref{ld} avoids the need to determine conservative contact forces, which, in the absence of dissipative contact forces, removes the need for free-body analysis.
For the case of a single rigid body, however, \eqref{ld} offers no advantage relative to a Newtonian-based derivation of Euler's equation.
Nevertheless, as an alternative derivation of \eqref{ee}, it is of interest to apply \eqref{ld} to a single rigid body.

A Lagrangian-based derivation of Euler's equation is given in \cite[p 281]{taeyoung} using Lagrangian dynamics on Lie groups.
The goal of the present note is to provide an elementary derivation by using generalized coordinates.
In particular, Euler angles are considered. 
Among all possible sequences consisting of three Euler-angle rotations, there are six that have three distinct axes and 6 that have the same first and last axes, for a total of 12 distinct sequences \cite[p. 764]{wertz}.
Relabeling axes allows us to consider two representative sequences, namely, 3-2-1 (azimuth-elevation-bank) and 3-1-3 (precession-nutation-spin).
These choices are commonly used for aircraft and spacecraft, respectively.
As a further illustration, Euler parameters (quaternions) are also considered.

The novelty of the present article is a compact derivation of Euler's equation using generalized coordinates.
Some elements of the derivation are known; for example, a connection is made between Proposition 1 and equation (A24) of \cite{meirovitch}.
However, a comparable derivation of Euler's equation using generalized coordinates does not appear to be available.


Notation:  For $x,y\in\BBR^3,$ $x\times y$ denotes the cross product of $x$ and $y,$ $x^\times$ denotes the cross-product matrix (so that $x^\times y = x\times y$), $I_3$ denotes the $3\times 3$ identity matrix, and $A^\rmT\in\BBR^{l\times k}$ denotes the transpose of $A\in\BBR^{k\times l}.$

\section{\large Preliminary Results}

For a single rigid body, let $q = [q_1\ q_2\ q_3]^\rmT\in\BBR^3$ denote generalized coordinates, and assume that the angular velocity $\omega$ can be parameterized as
\begin{equation}
    \omega(q,\dot q) = S(q)\dot q, \label{omegaS}
\end{equation}
where $S(q)\in\BBR^{3\times 3}.$
Assuming that the net force is zero and thus the center of mass of the spacecraft has zero inertial acceleration, it follows that
\begin{align}
    T(q,\dot q) &= \half \omega(q,\dot q)^\rmT J\omega(q,\dot q)\nn\\
    &= \half \dot q^\rmT S(q)^\rmT J S(q) \dot q,
\end{align}
and thus 
\begin{gather}
    \p_{\dot q}T(q,\dot q) = S(q)^\rmT J S(q)\dot q, \label{pdotT}\\
    \rmd_t\p_{\dot q}T(q,\dot q) = S(q)^\rmT J S(q)\ddot q
     + S(q)^\rmT J \dot S(q)\dot q + \dot S(q)^\rmT J S(q)\dot q, \label{dtpdotT}\\
    \p_{q}T(q,\dot q) 
    = \matl \dot q^\rmT [\p_{q_1}S(q)]^\rmT J S(q)\dot q\\
    \dot q^\rmT [\p_{q_2}S(q)]^\rmT J S(q)\dot q\\
    \dot q^\rmT [\p_{q_3}S(q)]^\rmT J S(q)\dot q\matr \label{pqdotT}.
\end{gather}
Furthermore, it follows from \cite[8.10.6]{baruh} that
\begin{equation}
    Q = S(q)^\rmT \tau. \label{QStau}
\end{equation}
Now, combining \eqref{dtpdotT}, \eqref{pqdotT}, \eqref{QStau} with \eqref{ld} yields
\begin{align}
    S(q)^\rmT J S(q)\ddot q
     + S(q)^\rmT J \dot S(q)\dot q + \dot S(q)^\rmT J S(q)\dot q 
    -  \matl \dot q^\rmT [\p_{q_1}S(q)]^\rmT J S(q)\dot q\\
    \dot q^\rmT [\p_{q_2}S(q)]^\rmT J S(q)\dot q\\
    \dot q^\rmT [\p_{q_3}S(q)]^\rmT J S(q)\dot q\matr = S(q)^\rmT \tau. \label{ld2}
\end{align}
If $S(q)$ is nonsingular, then
\begin{align}
    J S(q)\ddot q
     + J \dot S(q)\dot q +  S(q)^{-\rmT}\dot S(q)^\rmT J S(q)\dot q  
    -  S(q)^{-\rmT}\matl \dot q^\rmT [\p_{q_1}S(q)]^\rmT J S(q)]\dot q\\
    \dot q^\rmT [\p_{q_2}S(q)]^\rmT J S(q)\dot q\\
    \dot q^\rmT [\p_{q_3}S(q)]^\rmT J S(q)\dot q\matr =  \tau, \label{ld4}
\end{align}
which can be viewed as Euler's equation expressed in terms of arbitrary generalized coordinates.
Next, noting that
\begin{equation}
    \dot \omega(q,\dot q) = S(q)\ddot q + \dot S(q)\dot q,
\end{equation}
\eqref{ld4} can be written as
\begin{align}
    J \dot \omega(q,\dot q)
    +  S(q)^{-\rmT}\left(\dot S(q)^\rmT J S(q)\dot q  
    -   \matl \dot q^\rmT [\p_{q_1}S(q)]^\rmT J S(q)]\dot q\\
    \dot q^\rmT [\p_{q_2}S(q)]^\rmT J S(q)\dot q\\
    \dot q^\rmT [\p_{q_3}S(q)]^\rmT J S(q)\dot q\matr\right) =  \tau. \label{ld5}
\end{align}
Comparing \eqref{ld5} with Euler's equation \eqref{ee} written in terms of the angular velocity implies 
\begin{equation}
        S(q)^{-\rmT}\left(\dot S(q)^\rmT J S(q)\dot q  
    -    \matl \dot q^\rmT [\p_{q_1}S(q)]^\rmT J S(q)\dot q\\
    \dot q^\rmT [\p_{q_2}S(q)]^\rmT J S(q)\dot q\\
    \dot q^\rmT [\p_{q_3}S(q)]^\rmT J S(q)\dot q\matr\right)
    = \omega(q,\dot q)\times J\omega(q,\dot q). \label{eqnwanted}
    \end{equation}
Our objective is to verify this identity for rotations parameterized by Euler angles and Euler parameters.

For the following result, the columns of $S(q)$ are denoted by $S_1(q),$ $S_2(q),$ and $S_3(q)$ so that 
\begin{equation}
    S(q) = [S_1(q) \ \ S_2(q) \ \ S_3(q)].
\end{equation}

%
%
%
     %
    %
    %
    %
    %
%
%
%
    %
    %
    %
    %
%

We note that {\bf(a)} is given by equation (A24) of \cite{meirovitch}.

{\bf Proposition 1.}  Define $S$ by \eqref{omegaS}. Then, the following properties are equivalent:

{\bf(a)} For all $q$ and $\dot{q}$,
   \begin{equation}
	\dot S(q) + [S(q)\dot q]^\times S(q)  
	=    
	[ \p_{q_1}S(q)\dot q   \ \ \
	\p_{q_2}S(q)\dot q   \ \ \
	\p_{q_3}S(q)\dot q   ]. \label{identJ}
\end{equation}

{\bf(b)}  For all $q$ and $\dot{q}$,
\begin{align}
	\sum_{i=1}^3 \dot q_i\p_{q_i} S(q)   
	&+ [ S_2(q)\times S_3(q) \ \ \  S_3(q)\times S_1(q) \ \ \  S_1(q)\times S_2(q) ] \dot q^\times\nn\\
	&= [ \p_{q_1} S(q)\dot q \ \ \ \p_{q_2} S(q)\dot q \ \ \ \p_{q_3} S(q)\dot q]. 
	\label{identnoJ}
\end{align}

{\bf (c)} For all $q$,
 \begin{align}	 
 	\p_{q_2}S_1(q)-\p_{q_1}S_2(q)&=S_1(q)\times S_2(q),\label{identnodotq1}\\
 	\p_{q_3}S_1(q)-\p_{q_1}S_3(q)&=S_1(q)\times S_3(q),\label{identnodotq2}\\
 	\p_{q_3}S_2(q)-\p_{q_2}S_3(q)&=S_2(q)\times S_3(q).\label{identnodotq3} 
 \end{align}

Now, assume that $S(q)$ is nonsingular.  Then (a)--(c) are equivalent to
\begin{equation}\label{identshort}
S(q)^\rmT	[\p_{q_3}S_2(q)-\p_{q_2}S_3(q)\ \ \
	\p_{q_1}S_3(q)-\p_{q_3}S_1(q) \ \ \ 
	\p_{q_2}S_1(q)-\p_{q_1}S_2(q)]=\det(S(q))I_3. 
\end{equation}

The following lemmas are needed.

{\bf Lemma 1.} Let $x\in\BBR^3$ and $A\in\BBR^{3 \times 3}.$
Then,
\begin{align}
    A^\rmT (Ax)^\times A = (\det A)x^\times. \label{ATAxA}
\end{align}

{\bf Proof.}  See {\it xxxix}) of Fact 4.12.1 in \cite[p. 385]{bernstein2018}.
\hfill$\square$

{\bf Lemma 2.}  Let $A = [A_1\ \ A_2 \ \ A_3]\in\BBR^{3 \times 3}.$
Then
\begin{equation}
    A^\rmT[A_2\times A_3\ \ \ A_3\times A_1\ \ \ A_1\times A_2] 
    = (\det A)I_3. \label{ATA2I3}
\end{equation}
Now, let $x\in\BBR^3.$  Then, 
\begin{equation}
     [A_2\times A_3\ \ \ A_3\times A_1\ \ \ A_1\times A_2]x^\times 
    =   (Ax)^\times A. \label{A2A3AxA}
\end{equation}

{\bf Proof.}  See {\it xlii}) of Fact 4.12.1 in \cite[p. 385]{bernstein2018}.
In the case where $A$ is nonsingular, the second statement follows from \eqref{ATAxA} and \eqref{ATA2I3}.
In the case where $A$ is singular, the conclusion follows by continuity since both sides of \eqref{A2A3AxA} are continuous functions of the columns $(A_1,A_2,A_3)$ of $A$ and the set of nonsingular matrices is dense in $\BBR^{3\times3}$.\hfill$\square$

	

{\bf Proof of Proposition 1.}
Note that
\begin{equation}
    \dot S(q) = \sum_{i=1}^3\dot q_i\p_{q_i} S(q). \label{pf1}
\end{equation}
Furthermore, it follows from \eqref{A2A3AxA} that
\begin{align}
[S(q)\dot q]^\times S(q)   
&=   [ S_2(q)\times S_3(q) \ \ \  S_3(q)\times S_1(q) \ \ \  S_1(q)\times S_2(q) ]  \dot q^\times.      \label{pf2}
\end{align}
Therefore, \eqref{pf1} and \eqref{pf2} imply that {\bf(a)} and {\bf(b)} are equivalent.

To prove that  {\bf(b)} and {\bf(c)} are equivalent, note that {\bf(b)} is equivalent to $L(\dot{q})=R(\dot{q})$ for all $\dot{q}\in\BBR^3$, where $L$ and $R$ are the linear operators defined  for all $x=[x_1\ \ x_2\ \ x_3]^\rmT\in\BBR^3$ by
\begin{align}
L(x)&=
[ \p_{q_1} S(q)x \ \ \ \p_{q_2} S(q)x \ \ \ \p_{q_3} S(q)x]-	\sum_{i=1}^3  x_i\p_{q_i} S(q), \\ 
R(x)&= [ S_2(q)\times S_3(q) \ \ \  S_3(q)\times S_1(q) \ \ \  S_1(q)\times S_2(q) ] x^\times.
\end{align}
Since $R$ and $L$ are linear, it follows that $L(\dot{q})=R(\dot{q})$ for all $\dot{q}\in\BBR^3$ if and only if
\begin{equation}
L(e_i)=R(e_i),\quad\textrm{ for all $i=1,2,3$,}\label{eqLR}
\end{equation}
where
${e_1=[1\ \ 0\ \ 0]^\rmT}$,
${e_2=[0\ \ 1\ \ 0]^\rmT}$, and
${e_3=[0\ \ 0\ \ 1]^\rmT}$ because $(e_1,e_2,e_3)$ is a basis of $\BBR^3$.

Next, note that
 \begin{align}
 	L(e_1)&=  	[ \p_{q_1} S_1(q) \ \ \p_{q_2} S_1(q) \ \ \p_{q_3} S_1(q)]-
 	[ \p_{q_1} S_1(q) \ \  \p_{q_1} S_2(q) \ \ \p_{q_1} S_3(q)]\nn\\
 	&= [0 \ \ \ \p_{q_2} S_1(q)-\p_{q_1} S_2(q) \ \ \ \p_{q_3} S_1(q)-\p_{q_1} S_3(q)],\label{L1}\\
 	L(e_2)&= [\p_{q_1} S_2(q)-\p_{q_2} S_1(q) \ \ \ 0 \ \ \ \p_{q_3} S_2(q)-\p_{q_2} S_3(q)],\label{L2}\\
 	L(e_3)&= [\p_{q_1} S_3(q)-\p_{q_3} S_1(q) \ \ \ \p_{q_2} S_3(q)-\p_{q_3} S_2(q) \ \ \ 0],\label{L3}
 \end{align}
and
 \begin{align}
 R(e_1)&=[ S_2(q)\times S_3(q) \ \ \  S_3(q)\times S_1(q) \ \ \  S_1(q)\times S_2(q) ][0\ \ \ e_3 \ \ -e_2],\nn\\
 &=[0\ \  S_1(q)\times S_2(q) \ \ S_1(q)\times S_3(q)],\label{R1}\\
 R(e_2)&=[S_2(q)\times S_1(q) \ \ 0\ \ S_2(q)\times S_3(q)],\label{R2}\\
 R(e_3)&=[S_3(q)\times S_1(q) \ \ \ S_3(q)\times S_2(q)\ \ \ 0].\label{R3}
 \end{align}
Comparing \eqref{L1}--\eqref{L3} with \eqref{R1}--\eqref{R2} shows that \eqref{eqLR} is equivalent to {\bf(c)}. Finally \eqref{identshort} follows from \eqref{identnodotq1}--\eqref{identnodotq3} and \eqref{ATAxA}.
\hfill $\square$

To demonstrate the relevance of \eqref{identJ} to \eqref{eqnwanted}, note that transposing and rearranging \eqref{identJ} yields
\begin{equation}
     \dot S(q)^\rmT -    \matl \dot q^\rmT [\p_{q_1}S(q)]^\rmT   \\
     \dot q^\rmT [\p_{q_2}S(q)]^\rmT    \\
     \dot q^\rmT [\p_{q_3}S(q)]^\rmT    \matr
     =  S(q)^\rmT[S(q)\dot q]^\times,  
        \label{identJ3}
    \end{equation}
and thus, assuming that $S(q)$ is nonsingular,
\begin{equation}
     S(q)^{-\rmT}\left(\dot S(q)^\rmT -    \matl \dot q^\rmT [\p_{q_1}S(q)]^\rmT   \\
     \dot q^\rmT [\p_{q_2}S(q)]^\rmT    \\
     \dot q^\rmT [\p_{q_3}S(q)]^\rmT    \matr\right)
     =  [S(q)\dot q]^\times.  
        \label{identJ4}
    \end{equation}
Finally, multiplying \eqref{identJ4} on the right by $JS(q)\dot q$ yields
\begin{equation}
        S(q)^{-\rmT}\left(\dot S(q)^\rmT J S(q)\dot q  
    -    \matl \dot q^\rmT [\p_{q_1}S(q)]^\rmT J S(q)\dot q\\
    \dot q^\rmT [\p_{q_2}S(q)]^\rmT J S(q)\dot q\\
    \dot q^\rmT [\p_{q_3}S(q)]^\rmT J S(q)\dot q\matr\right)
    = \omega(q,\dot q)\times J\omega(q,\dot q),
    \end{equation}
which is precisely \eqref{eqnwanted}.

For a given choice of $q$, it easier to verify \eqref{identnodotq1}--\eqref{identnodotq3} than \eqref{identJ} or \eqref{identnoJ}.
In the next three sections, \eqref{identnodotq1}--\eqref{identnodotq3} are verified for 3-2-1 and 3-1-3 Euler angles as well as Euler parameters (quaternions). 

\section{\large Verification of \eqref{identnodotq1}--\eqref{identnodotq3} for 3-2-1 Euler Angles}

Letting $(\Psi,\Theta,\Phi)$ denote 3-2-1 (azimuth-elevation-bank) Euler angles, it follows that
\begin{equation}
\omega(q,\dot q)     = S(\Phi,\Theta)\dot q, \label{omegavectrixres}
\end{equation}
where
\begin{align}
S(\Phi,\Theta)
&=[  S_1 \ \ \  S_2(\Phi) \ \ \  S_3(\Phi,\Theta)]\\  
&=\matl 
1&0&-\sin\Theta \\
0&\cos\Phi &(\sin\Phi)\cos\Theta\\
0&-\sin\Phi &(\cos\Phi)\cos\Theta \matr, \label{Scols}
\end{align}
\begin{equation}
q =\matl q_1 \\     q_2 \\     q_3\matr\triangleq \matl \Phi \\     \Theta \\     \Psi\matr.
\end{equation}
Note that $\det S(\Phi,\Theta) = \cos\Theta,$ and thus $S(\Phi,\Theta)$ is singular if and only if gimbal lock occurs.
Hence,
%
%
\begin{align}	
	\p_{q_2}S_1(q)-\p_{q_1}S_2(q)&=	-\p_{\Phi}S_2(\Phi)=\matl 0\\ \sin\Phi \\ \cos\Phi\matr=	S_1\times S_2(\Phi),\nn\\
	\p_{q_3}S_1(q)-\p_{q_1}S_3(q)&=-\p_{\Phi}S_3(\Phi,\Theta)
	=\matl 0\\ -\cos\Phi\,\cos\Theta\\ \sin\Phi \,\cos\Theta\matr=
	S_1\times S_3(\Phi,\Theta),\\
	\p_{q_3}S_2(q)-\p_{q_2}S_3(q)&=-\p_{\Theta}S_3(\Phi,\Theta)=
	\matl \cos\Theta\\ \sin\Phi\,\sin\Theta\\ \cos\Phi \,\sin\Theta\matr=
		S_2(\Phi)\times S_3(\Phi,\Theta).\nn 
\end{align}
Hence, \eqref{identnodotq1}--\eqref{identnodotq3} hold, and thus \eqref{identJ}  and \eqref{identnoJ} are verified.

\section{\large Verification of \eqref{identnodotq1}--\eqref{identnodotq3} for 3-1-3 Euler Angles}

Letting $(\Phi,\Theta,\Psi)$ denote 3-1-3 (precession-nutation-spin) Euler angles, it follows that
\begin{equation}
\omega(q,\dot q)     = S(\Psi,\Theta)\dot q, \label{omegavectrixres313}
\end{equation}
where
\begin{align}
S(\Psi,\Theta)
&= [ S_1 \ \ \ S_2(\Psi) \ \ \ S_3(\Psi,\Theta)]\\
&= \matl 
0&\cos\Psi & (\sin\Psi) \sin\Theta\\
0&-\sin\Psi &(\cos\Psi)\sin\Theta  \\
1&0&\cos\Theta
\matr,
\end{align}
\begin{equation}
q=\matl q_1\\ q_2\\ q_3\matr\triangleq\matl 
   \Psi \\
   \Theta \\
   \Phi
\matr.
\end{equation}
Note that $\det S(\Psi,\Theta) = \sin\Theta,$ and thus $S(\Psi,\Theta)$ is singular if and only if gimbal lock occurs.
Hence,
%
\begin{align}	
	\p_{q_2}S_1(q)-\p_{q_1}S_2(q)&=	-\p_{\Psi}S_2(\Psi)=\matl \sin \Psi\\ \cos\Psi \\ 0\matr=	S_1\times S_2(\Psi),\nn\\
	\p_{q_3}S_1(q)-\p_{q_1}S_3(q)&=-\p_{\Psi}S_3(\Psi,\Theta)
	=\matl -\cos\Psi\,\sin \Theta,\\ \sin\Psi \,\sin\Theta \\ 0\matr=
	S_1\times S_3(\Phi,\Theta),\\
	\p_{q_3}S_2(q)-\p_{q_2}S_3(q)&=-\p_{\Theta}S_3(\Psi,\Theta)=
	\matl -\sin\Psi\,\cos\Theta\\ -\cos\Psi\,\cos\Theta\\ \sin\Theta\matr=
	S_2(\Psi)\times S_3(\Psi,\Theta).\nn 
\end{align}
Hence, \eqref{identnodotq1}--\eqref{identnodotq3} hold, and thus \eqref{identJ}  and \eqref{identnoJ} are verified.

 
\section{Verification of \eqref{identshort} for Euler Parameters}
To avoid gimbal lack, an alternative approach is to use Euler parameters (quaternions).
In this case,  
\begin{align}
\tilde q = \matl  q_1\\ q_2\\ q_3\\ q_4\matr =     \matl   \cos \half\theta\\[1ex]  (\sin\half\theta)n\matr, \label{eulerparvecgencoord}
\end{align}
where $\theta\in(-\pi,\pi]$ is the eigenangle and $n\in\BBR^3$ is the unit eigenaxis.
Since $q_1^2+q_2^2+q_3^2+q_4^2=1,$ it follows that $q_1 = \sqrt{1 - q_2^2-q_3^2-q_4^2},$ and thus the generalized coordinates are
$q = [q_2\ \ q_3\ \ q_4]^\rmT.$
With this notation, assuming that $\theta\ne\pi$ and thus $q_1>0,$ it follows that \eqref{omegaS} holds with
\begin{align}
S(q) =    2\matl 
q_1+\dfrac{q_2^2}{q_1} &
q_4+\dfrac{q_2q_3}{q_1} &
-q_3+\dfrac{q_2q_4}{q_1}\\[2ex]
-q_4+\dfrac{q_2q_3}{q_1} &
q_1+\dfrac{q_3^2}{q_1} &
q_2+\dfrac{q_3q_4}{q_1}\\[2ex]
q_3+\dfrac{q_2q_4}{q_1} &
-q_2+\dfrac{q_3q_4}{q_1} &
q_1+\dfrac{q_4^2}{q_1}\matr.
\end{align}
%

Next, note that, for all $i=2,3,4,$  $\p_{q_i}q_1=-q_i/q_1$. Thus,
\begin{align}
	\p_{q_3}S_1(q)-\p_{q_2}S_2(q)&=2\matl
	-\dfrac{q_3}{q_1}+\dfrac{q_2^2q_3}{q_1^3}-\dfrac{q_3}{q_1}-\dfrac{q_2^2q_3}{q_1^3}\\[2ex]
	\dfrac{q_2}{q_1}+\dfrac{q_2q_3^2}{q_1^3}+\dfrac{q_2}{q_1}-\dfrac{q_3^2q_2}{q_1^3}\\[2ex]
	1+\dfrac{q_2q_3q_4}{q_1^3}+1-\dfrac{q_2q_3q_4}{q_1^3}
	\matr=\frac{4}{q_1}\matl
	-q_3\\
	q_2\\
    q_1
	\matr,\\
	\p_{q_4}S_1(q)-\p_{q_2}S_3(q)&=2\matl -\dfrac{q_4}{q_1}+\dfrac{q_2^2q_4}{q_1^3}-\dfrac{q_4}{q_1}-\dfrac{q_2^2q_4}{q_1^3}\\[2ex]
	-1+\dfrac{q_2q_3q_4}{q_1^3}-1-\dfrac{q_2q_3q_4}{q_1^3}\\[2ex]
	\dfrac{q_2}{q_1}+\dfrac{q_2q_4^2}{q_1^3}+\dfrac{q_2}{q_1}-\dfrac{q_2q_4^2}{q_1^3}
	\matr=\frac{4}{q_1}\matl
	-q_4\\
	-q_1\\
 	 q_2
	\matr,\\
	\p_{q_4}S_2(q)-\p_{q_3}S_3(q)&=2\matl
	1+\dfrac{q_2q_3q_4}{q_1^3}+1-\dfrac{q_2q_3q_4}{q_1^3}\\[2ex]
	-\dfrac{q_4}{q_1}+\dfrac{q_3^2q_4}{q_1^3}-\dfrac{q_4}{q_1}-\dfrac{q_3^2q_4}{q_1^3}\\[2ex]
	\dfrac{q_3}{q_1}+\dfrac{q_3q_4^2}{q_1^3}+\dfrac{q_3}{q_1}-\dfrac{q_3q_4^2}{q_1^3}
	\matr=\frac{4}{q_1}\matl
	q_1\\
	-q_4\\
	q_3
	\matr.
\end{align}
Thus, 
\begin{align}
	M(q)&\triangleq[\p_{q_4}S_2(q)-\p_{q_3}S_3(q)\ \ 
		\p_{q_2}S_3(q)-\p_{q_4}S_1(q)\ \ 
	\p_{q_3}S_1(q)-\p_{q_2}S_2(q)]\nn\\
	&=
	\frac{4}{q_1}\matl
 q_1&q_4&-q_3\\ -q_4&q_1&q_2 \\q_3&-q_2&q_1	\matr,
\end{align}
and using $q_1^2+q_2^2+q_3^2+q_4^2=1$ yields
\begin{align}
	M(q)^\rmT S(q)&=
	\frac{8}{q_1}\matl
	q_1&-q_4&q_3\\ q_4&q_1&-q_2 \\-q_3&q_2&q_1	\matr
\matl 
q_1+\dfrac{q_2^2}{q_1} &
q_4+\dfrac{q_2q_3}{q_1} &
-q_3+\dfrac{q_2q_4}{q_1}\\[2ex]
-q_4+\dfrac{q_2q_3}{q_1} &
q_1+\dfrac{q_3^2}{q_1} &
q_2+\dfrac{q_3q_4}{q_1}\\[2ex]
q_3+\dfrac{q_2q_4}{q_1} &
-q_2+\dfrac{q_3q_4}{q_1} &
q_1+\dfrac{q_4^2}{q_1}\matr	\nn\\
&=\frac{8}{q_1}\matl 1&0&0\\0&1&0\\0&0&1\matr. \label{eqMS}
\end{align}
Since $\det(M(q))=64/q_1^2$, \eqref{eqMS} implies that $\det(S(q))=8/q_1$. Thus, $S(q)$ is nonsingular and satisfies \eqref{identshort}.
Consequently, \eqref{identJ}  and \eqref{identnoJ} are verified for Euler parameters.

\section*{Acknowledgments}
The authors are grateful to one of the reviewers for bringing \cite{meirovitch} to our attention and independently confirming {\bf (c)} of Proposition 1.

\section{Conclusions}

We used Lagrangian dynamics to derive Euler's equation using Euler angles and Euler parameters (quaternions) as generalized coordinates.
Although the strength of Lagrangian dynamics lies in its ability to avoid free-body analysis in the presence of conservative reaction forces, this derivation illustrates the connection between Lagrangian dynamics and the dynamics of a single unconstrained rigid body.
A more advanced approach is to apply Lagrangian dynamics on Lie groups as presented in \cite{taeyoung}.

\bibliographystyle{myIEEEtran}
\bibliography{euler}

\end{document}